\documentclass[12pt]{amsart}
\usepackage{amscd,amssymb}
\usepackage[arrow,matrix]{xy}

\theoremstyle{plain}
\newtheorem{theorem}[subsection]{Theorem}
\newtheorem{lemma}[subsection]{Lemma}

\newtheorem{corollary}[subsection]{Corollary}

\theoremstyle{definition}
\newtheorem{remark}[subsection]{Remark}
\newtheorem{definition}[subsection]{Definition}
\newtheorem{example}[subsection]{Example}

\numberwithin{equation}{section}
\newcommand{\X}{{\mathcal X}}
\newcommand{\Y}{{\mathcal Y}}
\newcommand{\ZZ}{{\mathcal Z}}
\newcommand{\T}{{\mathcal T}}
\newcommand{\F}{{\mathcal F}}
\newcommand{\A}{{\mathcal A}}
\newcommand{\wA}{\widehat{{\mathcal A}}}
\newcommand{\B}{{\mathcal B}}
\newcommand{\LL}{{\mathcal L}}
\newcommand{\D}{{\mathcal D}}

\renewcommand{\SS}{{\mathcal S}}

\newcommand{\Z}{\mathbb{Z}}
\newcommand{\Q}{\mathbb{Q}}

\newcommand{\C}{\mathbb{C}}
\newcommand{\K}{\mathbb{K}}
\newcommand{\PP}{\mathbb{P}}
\newcommand{\bF}{\mathbb{F}}

\newcommand{\tX}{\widetilde{X}}
\newcommand{\tY}{\widetilde{Y}}
\newcommand{\tZ}{\widetilde{Z}}
\newcommand{\tM}{\widetilde{M}}

\DeclareMathOperator{\rank}{rank}

\DeclareMathOperator{\coker}{coker}

\DeclareMathOperator{\id}{id}

\DeclareMathOperator{\tor}{Tor}

\newcommand{\surj}{\twoheadrightarrow}

\begin{document}
\date{May 14, 2003}

\title[Equivariant chain complexes of arrangements]
{Equivariant chain complexes, twisted homology\\
and relative minimality of arrangements}

\author[Alexandru Dimca]{Alexandru Dimca}
\address{Laboratoire de Math. Pures, Univ.
Bordeaux I, 33405 Talence Cedex, France}
\email
{dimca@math.u-bordeaux.fr}

\author[\c{S}tefan Papadima]{\c{S}tefan Papadima$^*$}
\address{Inst. of Math. "Simion Stoilow",
P.O. Box 1-764,
RO-70700 Bucharest, Romania}
\email
{Stefan.Papadima@imar.ro}

\thanks{$^*$Partially supported by grant 
CNCSIS 693/2002 of the Romanian Ministry of
Education and Research}



\begin{abstract}
We show that the $\pi$-equivariant chain complex
($\pi=\pi_1(M(\A))$),
$C_{\bullet}(\tX)$, associated to a Morse-theoretic
minimal $CW$-structure $X$ on the complement $M(\A)$ of
an arrangement $\A$, is independent of $X$. The same holds
for all scalar extensions, $C_{\bullet}(\tX)\otimes_{\Z \pi}
\K \Z$, $\K$ a field, where $X$ is an arbitrary minimal
$CW$-structure on a space $M$. When $\A$ is a section of another arrangement
$\wA$, we show that the divisibility properties of the first
Betti number of the Milnor fiber of $\A$ obstruct the 
homotopy realization of $M(\A)$ as a subcomplex of a
minimal structure on $M(\wA)$.

If $\wA$ is aspherical and $\A$ is a sufficiently 
generic section of $\wA$, then $H_*(M(\A) ;L)$ 
may be described in terms of $\pi$, $L$ and $\chi (M(\A))$,
for an arbitrary local system $L$; explicit computations 
may be done, when $\wA$ is fiber-type. In this case, 
explicit $\K \Z$-presentations of arbitrary abelian
scalar extensions
of the first non-trivial higher homotopy group of $M(\A)$,
$\pi_p$, may also be obtained. For nonresonant abelian
scalar extensions,
the $\C \Z$-rank of $\pi_p\otimes_{\Z \pi} \C \Z$ 
is combinatorially determined.
\end{abstract}

\maketitle

\tableofcontents

\section{Introduction}
\label{sec:intro}

Let $\A =\{ H_0, H_1, \dots, H_n \}$ be a complex 
hyperplane arrangement in $\PP^{r-1}$, with complement
$M=M(\A)=\PP^{r-1}\setminus \cup_{i=0}^n H_i$, and
fundamental group $\pi =\pi_1(M)$. The cohomology 
ring of the complement,
with arbitrary {\em constant} coefficients, was computed 
by Orlik-Solomon~\cite{OS}; their description involves
solely the {\em combinatorics} of $\A$, that is, the
associated intersection lattice $\LL(\A)$. The fundamental group
$\pi$ is complicated in general, but nevertheless rather well
understood and accessible to concrete computations, see for
instance \cite{CS3}, \cite{L1}. Note however that $\LL(\A)$
does {\em not} determine $\pi$ in general, as the example of
Rybnikov~\cite{Ry} shows.

Much less is known about the homology groups $H_*(M;L)$ of $M$ with 
{\em twisted} coefficients (alias $\Z \pi$--modules), $L$.
Our aim in this paper is to get more insight on such groups
$H_*(M;L)$, both at the general and computational 
level, along the lines sketched in our previous 
work~\cite[Remark 12(ii)]{DP}.

\subsection{Absolute minimality and equivariant chain complexes}
\label{i1=subsec}

Among other things, we have proved in~\cite{DP} that $M$
has the homotopy type of a {\em minimal} $CW$-complex $X$,
that is, of a $CW$-complex with trivial cellular incidences
(a result independently obtained by Randell~\cite{R}).
In this note, we go further and analyze in Section~\ref{sec:mineq}
the associated $\Z \pi$--chain complex of the universal cover
$C_{\bullet}(\tX)$. As it is well-known~\cite{W}, this is the
{\em universal} object encoding twisted homology information,
since
\begin{equation}
\label{univ=eq}
H_*(M;L)=H_*(C_{\bullet}(\tX)\otimes_{\Z \pi} L), \quad 
\text{for all} \quad L \, .
\end{equation}
At the same time, this is also relevant for computations
related to {\em higher homotopy groups} 
(a very difficult subject, in general), since
\begin{equation}
\label{hipi=eq}
\pi_p(M)=H_p(C_{\bullet}(\tX))\, ,
\end{equation}
by Hurewicz, where $\pi_p(M)$ is the first higher non-trivial
homotopy group of $M$ (for $M$ non-aspherical).

Note that the $\Z \pi$--chain complex $C_{\bullet}(\tX)$
may depend on the choice of $CW$-structure on $M$. 
In general, all one can say is that the 
$\Z \pi$--{\em chain homotopy type}
is independent of $X$. In this direction, 
we have two general results where
we can upgrade the above independence property to
{\em chain isomorphism type}.

The first one is Corollary~\ref{cor:indarr}, which says
that the $\Z \pi$--chain complex $C_{\bullet}(\tX)$
actually defines a {\em new} invariant of $\A$, depending
only on the projective equivalence type of $\A$. Key to
our proof is an argument involving Whitney stratifications
and transversality conditions for projective flags. In this
way, we recover in particular a series of results, proved by
Cohen--Orlik in~\cite{C} and~\cite{CO}.

A standard technique in topology is  abelianization
(leading for instance to the theory of Alexander invariants, see~\cite{Hi}).
When a $\Z \pi$--module like $\pi_p(M)$ in~\eqref{hipi=eq} above
looks intractable, it is generally hoped that
$\pi_p(M)\otimes_{\Z \pi}\K \Z$ will be more manageable, where
$\K$ is a field and
\begin{equation}
\label{ab=eq}
\nu : \pi \longrightarrow \Z
\end{equation}
is a group character, inducing a change of rings,
$\nu : \Z \pi \to \K \Z$. The abelian extensions of scalars 
$\pi_p(M)\otimes_{\Z \pi}\K \Z$  obtained in this way will be
called {\em character--abelianizations}.

Our second main result in Section~\ref{sec:mineq},
Theorem~\ref{thm:indmin}, says that the $\K \Z$--chain complex
$C_{\bullet}(\tX)\otimes_{\Z \pi} \K \Z$ is independent of the
{\em minimal} $CW$-structure $X$ on $M$, for an {\em arbitrary} 
space $M$.

\subsection{General position and cell structures}
\label{i2=subsec}

Now assume that $\A$ is a $k$--{\em generic section} of
another arrangement, $\wA$ in $\PP^{m-1}$. By `$k$--genericity'
we mean a set of transversality conditions, depending on $k$,
with respect to certain strata of the natural stratification of
$\PP^{m-1}$ associated to $\wA$; see~\cite[(1)]{DP} for details.
When we speak about $\wA$--{\em genericity}, we simply mean that
$\PP^{r-1}$ is transverse to all $\wA$--strata.

Following our program from~\cite[\S \S 5--6]{DP}, our idea here
is to get from $M(\wA)$ twisted homology information on
$M(\A)$, assuming $k$--genericity. The key tool is provided 
by~\cite[Proposition 14]{DP}, which says that one may replace,
up to homotopy, the inclusion $M(\A)\hookrightarrow M(\wA)$ by
a cellular map, $f: X \to Y$, between minimal $CW$-complexes,
such that
\begin{equation}
\label{kapp=eq}
f_{\mid X^{(k)}}= \id \, .
\end{equation}

When $M(\wA)$ is aspherical and $k\ge 2$, it follows 
from~\eqref{kapp=eq} above
that $M(\wA)$ is a $K(\pi, 1)$. Moreover, Proposition 15
from~\cite{DP} implies that one may replace in~\eqref{kapp=eq}
$k$ by $p= p(M(\A))$, where $p$ is a homotopy invariant of
$M(\A)$, introduced in~\cite{PS}.

Assuming full $\wA$--genericity, \eqref{kapp=eq} improves to
\begin{equation}
\label{skapp=eq}
X =Y^{(r-1)}\, ,
\end{equation}
where $f$ is the inclusion of the $(r-1)$--skeleton,
$Y^{(r-1)}\hookrightarrow Y$. The basic example when
the homotopy formula~\eqref{skapp=eq} holds is provided by
Boolean genericity (in other words, usual general position);
see Hattori's pioneering paper~\cite{H}. This generalizes to
{\em fiber-type genericity}. The arrangements $\A$ which are
$\wA$--generic sections of fiber-type arrangements $\wA$ may be
defined in purely combinatorial terms: they are the {\em hypersolvable}
arrangements (introduced in~\cite{JP1}) for which $p=r-1$; see
\S~\ref{subsec:genhat}.

\subsection{Relative minimality and twisted homology}
\label{i3=subsec}

In Section~\ref{sec:relmin} (Theorem~\ref{thm:best} and
Example~\ref{ex:no}) we prove that, in general, the 
homotopy formula~\eqref{kapp=eq} from \S~\ref{i2=subsec} is
the best possible result. In spite of some 
homological and Morse-theoretic    
positive evidence, it turns out that the ideal
{\em relative minimality} formula,
\begin{equation}
\label{ideal=eq}
f \quad = \quad \text{subcomplex inclusion}\, ,
\end{equation}
cannot hold in full generality, not even on $X^{(k+1)}$.

To prove this, we explore the consequences 
of~\eqref{ideal=eq} on the first Betti number of the
Milnor fiber of $\A$. This in turn uses a decomposition 
of the homology of the Milnor fiber, found by Cohen--Suciu
in~\cite{CS1}, in terms of the homology of $M(\A)$ with certain
twisted coefficients. 

\subsection{Fiber-type arrangements and free differential calculus}
\label{i4=subsec}

Most results on twisted homology of arrangements available in the
literature assume {\em nonresonant} coefficients and take the form of 
vanishing theorems. There is however a particularly remarkable 
exception, due to Cohen--Suciu~\cite{CS2}, from the computational
point of view.

The authors of~\cite{CS2} consider finitely-presented groups, $\pi$,
which are iterated semidirect products of free groups, with all
monodromy actions trivial on homology. They use the Fox
{\em free differential calculus} to construct a minimal
$CW$-structure $Y$ on $K(\pi, 1)$, and to describe {\em explicitly}
the associated $\Z \pi$--chain complex, $C_{\bullet}(\tY)$.

On one hand, this leads to explicit twisted homology computations
with non-trivial result; see~\cite{CS1},~\cite{CS2}. 
On the other hand, this has implications in arrangement theory,
coming from the fact that the fundamental groups 
of fiber-type arrangements
all have such a semidirect product structure; see Falk--Randell~\cite{FR}.

\subsection{Twisted homology with resonant coefficients}
\label{i5=subsec}

We come back, in Section~\ref{sec:tord}, to our framework,
explained in \S~\ref{i2=subsec}. Let $\A$ be a 
sufficiently generic section
of an aspherical arrangement $\wA$, with complement 
$M=M(\A)$ and fundamental group $\pi =\pi_1(M)$. Set
$p=p(M)$. Let $L$ be a {\em resonant} (that is, arbitrary)
$\Z \pi$--module.

Using~\eqref{kapp=eq}, we infer in Theorem~\ref{thm:range} that
\begin{equation}
\label{tor=eq}
H_{<p}(M;L) =\tor_{<p}^{\Z \pi}(\Z, L)\, .
\end{equation}
When $\wA$ is fiber-type, the $\tor^{\Z \pi}$--groups
from~\eqref{tor=eq} above may be explicitly computed by
Fox calculus; see \S~\ref{i4=subsec}. Assuming only that
there are no collinearity relations among the defining equations
of the hyperplanes of $\A$, Theorem~\ref{thm:range} applies
and everything becomes extremely simple:
$\pi =\Z^n$, and the Fox resolution is the standard Koszul
resolution of $\Z$ over $\Z \Z^n$; see Corollary~\ref{cor:ccomb}.

In Theorem~\ref{thm:compl}, we assume full $\wA$--genericity,
hence the stronger homotopy formula~\eqref{skapp=eq} becomes available.
Via an Euler characteristic argument, we are thus able to extend
\eqref{tor=eq} above to a full computation of $H_*(M;L)$, 
involving only $\tor_*^{\Z \pi}(\Z, L)$ and $\chi (M)$.

\subsection{Character--abelianizations of higher homotopy groups}
\label{i6=subsec}

Very little is known about higher homotopy groups of
complements in $\PP^{r-1}$ of complex projective hypersurfaces, 
$V(h) = \{ h=0 \}$. Let $H$ be a hyperplane in $\PP^{r-1}$,
generic with respect to $V(h)$. 
Set $M= \PP^{r-1}\setminus (V(h) \cup H)$, and $\pi =\pi_1(M)$.
When $h$ is irreducible,  Libgober~\cite{L2} showed that,
under certain additional hypotheses, $\pi =\Z$, and $\pi_p(M)$
(the first higher non-trivial homotopy group of $M$) is a
torsion $\Q \Z$--module.

Arrangement complements are the simplest {\em non}--irreducible
examples. Here, $\pi$ is much more complicated, and 
character--abelianizations of $\pi_p(M)$ may have strictly positive 
$\Q \Z$--rank; see our Theorem~\ref{thm:rcomb}~\eqref{d3}
from Section~\ref{sec:ab}.

Our results on character--abelianizations of $\pi_p$ assume full
$\wA$--genericity, like in Theorem~\ref{thm:compl}
(see \S~\ref{i5=subsec}). The reason is that in this case
the minimal $\Z \pi$--resolution of $\pi_p(M)$ from
\cite[Theorem 18(ii)]{DP} is available.

In Theorem~\ref{thm:rarb}, we give a formula for the
$\K \Z$--rank of an arbitrary 
character--abelianization of $\pi_p(M)$,
involving $\tor_*^{\Z \pi}(\Z, \K \Z)$ and $\chi (M)$.

Assume now moreover that $\wA$ is fiber-type. Then
Corollary~\ref{cor:piexpl} provides an explicit
$\K \Z$--{\em presentation} of 
an arbitrary character--abelianization,
$\pi_p(M)\otimes_{\Z \pi}\K \Z$. In this case,
Theorem~\ref{thm:rcomb} also provides an explicit 
{\em combinatorial} formula,
for the $\C \Z$--rank of a so-called $\wA$--nonresonant
character--abelianization (see Definition~\ref{def:nr}) of $\pi_p(M)$.
The last result represents a new combinatorial determination
phenomenon, related to higher homotopy groups of arrangements,
to be compared with~\cite[Theorem 23]{DP}.
 
\section{Minimal equivariant chain complexes}
\label{sec:mineq}

\subsection{Generic flags and minimal $CW$-structures}
\label{subsec:genflag}

Let $\A$ be a complex hyperplane arrangement in $\PP (U)$,
with complement $M:= M(\A)$ and fundamental group
$\pi:= \pi_1(M)$. In Section $4$ from \cite{DP}, we have
constructed a {\em minimal} $CW$-structure on $M$, that is,
a $CW$-complex $X$ homotopy equivalent to $M$, having as many 
$k$-cells as the $k$-th Betti number, for all $k$. Our
construction involved various choices. To examine this issue,
we will need several definitions. Without any loss of
generality, we may assume that $\A$ is {\em essential}; see
\cite[p.197]{OT}.

The {\em intersection lattice}, $\LL (\A)$, is the set of
{\em edges} of $\A$ (that is, the nonempty intersections of 
hyperplanes from $\A$), ordered by reverse inclusion.
One has a canonical stratification of   $\PP (U)$,
$\SS (\A):= \{ \SS_{S} \}_{S\in \LL (\A)}$, with the property that
$\overline{\SS_{S}}= S$, for all $S$; see \cite[III.3.1 and III.4.5]
{GM}. Set $r:= \dim \, U$.  

\begin{definition}
\label{def:flag}

An increasing projective flag in $\PP (U)$, 
$\F = \{ \F_s \}_{0\le s<r}$
(where $\dim \, \F_s =s$, for all $s$)
is called $\A$--{\em generic} if $\F_s$ is transverse
to $\SS (\A^{\F_{s+1}})$, for $s<r-1$, where
$\A^{\F_{s+1}}$ denotes the projective arrangement in 
$\F_{s+1}$ obtained by restriction from $\A$.
\end{definition}

An $\A$--generic flag gives a finite increasing
filtration of $M$, $\{ F_s M \}_{0\le s<r}$:
\begin{equation}
\label{eq:ffilt}
F_s M:= M \cap \F_s \, .
\end{equation}

\begin{definition}
\label{def:minad}

Let $X$ be an $(r-1)$--dimensional minimal complex.
A homotopy equivalence, $\varphi : X\to M$, is {\em compatible}
with a given $\A$--generic flag if $\varphi$ restricts to homotopy
equivalences, $\varphi_s : X^{(s)} \to F_s M$, for $0\le s<r$. 
\end{definition}

Our (absolute) minimality result from \cite[Section 4]{DP}
reads then: for any $\A$--generic flag, $\F$, there exists an
$\F$--compatible minimal structure on $M(\A)$, $X$, in the sense of
the above definition \ref{def:minad}.

Let now $(Z, z_0)$ be a reasonable based space (e.g., a 
connected complex, or a connected manifold), endowed with a finite 
increasing filtration, $\{ F_s Z \}_{s\ge 0}$, such that 
$F_0 Z= \{ z_0 \}$. Let $p_Z : (\tZ ,\tilde{z_0}) \to (Z, z_0)$
be a universal cover. Set $\pi:= \pi_1(Z, z_0)$. Define
\[
F_s \tZ := p_Z^{-1} (F_s Z)\, , \quad \text{for} \quad s\ge 0 \, .
\]
The associated $\pi$--{\em equivariant chain complex},
\begin{equation}
\label{eq:nottw}
C_{\bullet}(\tZ ):= \{ d_q : H_q(F_q \tZ ,F_{q-1}\tZ) 
\longrightarrow H_{q-1}(F_{q-1}\tZ ,F_{q-2}\tZ ) \}_q \, ,
\end{equation}
is endowed with the boundary maps, $d_q$, coming from 
the triple $(F_q, F_{q-1}, F_{q-2})$, in the standard way.
It is a chain complex of (right) $\Z \pi$--modules.

We are going to show that the (minimal) $\pi$--equivariant
chain complex, $C_{\bullet}(\tX)$, from \cite{DP}, is a 
well-defined invariant of $\A$. We do this in two steps.

Let $\varphi : X\to M$ be an $\F$--compatible homotopy
equivalence, as in Definition~\ref{def:minad}. Denote by
$\Phi$ the lift of $\varphi$ to (based) universal covers.
Use $\varphi_{\#} : \pi_1(X)\stackrel{\sim}{\to} \pi_1(M)$
to identify fundamental groups. Denote by
\begin{equation}
\label{eq:comp}
\widetilde{\varphi} := \{ \widetilde{\varphi_q}:
C_q (\tX)\longrightarrow C_q(\tM) \}_q
\end{equation}
the $\Z \pi$--chain map, induced by $\Phi$, between the
equivariant chain complexes associated to the skeletal filtration
on $X$, and the $\F$--filtration \eqref{eq:ffilt} on $M$
respectively.

\begin{lemma}
\label{lem:ind1}

The above map~\eqref{eq:comp} is an isomorphism of
$\Z \pi$--chain complexes.
\end{lemma}

\begin{proof}
By standard homotopy properties of fibrations,
$p_X : (\tX, F_q \tX)\to (X, F_q X)$ will induce isomorphisms
on all homotopy groups, for $q\ge 0$, and likewise for $M$.
It follows that $\Phi : (\tX, F_q \tX)\to (\tM, F_q \tM)$
(and therefore $\Phi : F_q \tX \to F_q \tM$) induces isomorphisms
on all homotopy groups, since $\varphi$ induces homotopy
equivalences, $\varphi_q : F_q X \to F_q M$. By Whitehead's
theorem, $\Phi : (F_q \tX, F_{q-1} \tX) \to (F_q \tM, F_{q-1} \tM)$
induces homology isomorphisms, for all $q$.
\end{proof}

This shows that the skeletal equivariant chain complex does not
depend on the choice of $\F$--compatible minimal structure, for 
a fixed $\A$--generic flag, $\F$.

\begin{lemma}
\label{lem:ind2}

Let $\F$ and $\F'$ be two $\A$--generic flags, with associated 
filtrations $F_{\bullet}M$ and $F_{\bullet}'M$ respectively.
There is a filtered homeomorphism,
\[
\Psi : (M, F_{\bullet} M) \stackrel{\sim}{\longrightarrow} 
(M, F_{\bullet}' M)\, .
\]
\end{lemma}

\begin{proof}
The set $Flag(\PP (U))$ of all the flags
$$ \F: \emptyset=\F_{-1} \subset \F_0 \subset ... \subset \F_{r-1}$$
in $\PP (U)$ 
is a complex algebraic variety which is compact, smooth and connected,
hence irreducible. The subset $Flag(\A) \subset Flag(\PP (U))$
of all the $\A$-generic flags is an open Zariski subset in this 
irreducible variety, and as such it is connected.
Consider the second projection
$$p_2:\PP (U) \times Flag(\A) \to Flag(\A).$$
On the product $\PP (U) \times Flag(\A)$ there are two natural Whitney regular stratifications. 

The first one, denoted by $\X$, is the product of the stratification $\SS(\A)$ on  $\PP (U)$ by the trivial stratification on $ Flag(\A)$, i.e. $\X$ has as strata the products of the form
$$X_S=\SS _S \times  Flag(\A)$$
for $S \in \LL(\A)$, see \cite{GWPL}, p.12.
 
The second stratification, denoted by $\Y$, has as strata the following (constructible) submanifolds in $\PP (U) \times Flag(\A)$:
$$Y_j=\{(x,\F) \in \PP (U) \times Flag(\A);x \in \F_j \setminus \F_{j-1} \}$$
for $j=0,1,...,r-1.$

To see that the  stratification $\Y$ is Whitney regular, note that it is induced by the stratification $\ZZ$ on $\PP (U) \times Flag(\PP (U))$ given by the orbits of the regular action
$$Aut(U) \times (\PP (U) \times Flag(\PP (U))) \to \PP (U) \times Flag(\PP (U))$$
where $g \cdot (x,\F)=(gx,g\F)$. Indeed, the orbits of this action are exactly the sets
$$Z_j=\{(x,\F) \in \PP (U) \times Flag(\PP (U)  );x \in \F_j \setminus \F_{j-1} \}$$
for $j=0,1,...,r-1.$
The resulting stratification $\ZZ$ is  Whitney regular by a general result
on regular actions with finitely many orbits, see \cite{GWPL}, p.21.

Using the definition of an $\A$-generic flag, it follows that the two stratifications $\X$ and $\Y$ are transversal (or in general position) and hence their intersection $\T$ is also a  Whitney regular stratification, see \cite{GWPL}, p.12. Note that the strata
of this stratification are all the non-empty intersections $T_{S,j}=X_S \cap Y_j$, in particular these intersections are connected. It can be shown by a
direct computation that the restrictions
$$p_2:T_{S,j} \to  Flag(\A)$$
are all submersive surjections. Applying Thom's First Isotopy Lemma,
 see Theorem (5.2), \cite{GWPL}, p.58, to the second projection $p_2$ we get that
$\PP (U) \times Flag(\A)$ fibers over $ Flag(\A)$ in the stratified sense.
This means that for any two fibers of $p_2$, corresponding to two $\A$-generic
flags $\F$ and $\F'$, there is a homeomorphism $\PP (U)\times \{\F\}  \to \PP (U) \times \{\F'\} $ sending
the strata of the stratifications induced by $\T$ into each other. This proves our claim by considering only the strata contained in $M$.
\end{proof}

\begin{corollary}
\label{cor:indarr}

The $\Z \pi$--chain complex associated to any Morse-theoretic
minimal $CW$-structure, $X$, on $M(\A)$, constructed in \cite{DP},
is an invariant of $\A$. Actually, the $\Z \pi$--chain complex 
$C_{\bullet}(\tX)$ depends only on the projective equivalence type
of $\A$.
\end{corollary}

\subsection{Relation to work by Cohen--Orlik} 
\label{subsec:co}

The $\pi$--equivariant chain complex from 
Corollary~\ref{cor:indarr}
is a potentially powerful invariant of the arrangement $\A$.
This is due for instance to the fact that
the first nonzero homology group of $C_{\bullet}(\tX)$ is 
isomorphic to
the first nonzero higher homotopy group of $M(\A)$, when 
$M(\A)$ is not aspherical; see~\cite[\S\S 5--6]{DP},
for various results on higher homotopy groups, inspired 
from this remark.
This is also due to the fact that  $C_{\bullet}(\tX)$ is
the {\em universal} chain complex
computing the twisted (co)homology of the complement. Indeed,
let $N$ be an arbitrary left $\Z \pi$--module. Then
\begin{equation}
\label{eq:spec}
H_*(M(\A); N)= H_*(C_{\bullet}(\tX)\otimes_{\Z \pi} N)
\end{equation}
(and similarly for cohomology); see \cite[Ch. VI]{W}.

From basic equation~\eqref{eq:spec}, we may easily recover
(in homological form) various results, found by D. Cohen in
\cite{C} for the case of a $\C$--vector space $N$, and 
then reformulated in terms of flags by D. Cohen and P. Orlik 
\cite{CO}. Note that the flags $\F$ used in \cite[\S 2.3]{CO} coincide
with our $\A$--generic flags from Definition~\ref{def:flag}.

If $(Z, z_0)$ is a reasonable filtered space, 
as in \S~\ref{subsec:genflag}, 
and $N$ is a left $\Z \pi_1(Z, z_0)$--module, one may consider 
the chain complex 
\begin{equation}
\label{eq:tw}
C_{\bullet}(Z ;N):= \{ d_q : H_q(F_q Z, F_{q-1}Z; N)
\longrightarrow H_{q-1}(F_{q-1}Z, F_{q-2}Z; N) \}_q \, ,
\end{equation}
(the twisted version of \eqref{eq:nottw}). If 
$\varphi: X \to M$  
is an $\F$--compatible homotopy equivalence, as in Definition
\ref{def:minad}, it induces a chain isomorphism,
\begin{equation}
\label{eq:mx}
C_{\bullet}(X; N)\stackrel{\sim}{\longrightarrow}
C_{\bullet}(M; N)\, ,
\end{equation}
where $C_{\bullet}(M; N)$ is the homology version of the
complex from \cite{C} and \cite{CO}. On the other hand,
\begin{equation}
\label{eq:spec2}
C_{\bullet}(X; N)\simeq C_{\bullet}(\tX)\otimes_{\Z \pi} N
\end{equation}
(isomorphism of chain complexes), for arbitrary $N$; see
\cite[Theorem VI.4.9]{W}. At the same time,
\begin{equation}
\label{eq:dim}
\dim_{\C} ( C_q (\tX)\otimes_{\Z \pi} N)=
(\dim_{\C} N)\cdot b_q(M(\A))\, ,\quad \forall q 
\end{equation}
(from minimality), if $N$ is a finite-dimensional $\C$-vector space.

Equations~\eqref{eq:spec}, \eqref{eq:spec2} and~\eqref{eq:mx}
together imply that the chain complex $C_{\bullet}(M; N)$ 
computes the twisted homology of the complement, $H_*(M(\A); N)$;
if $\dim_{\C} N=1$, \eqref{eq:dim} implies that, furthermore,
$\dim_{\C} C_q(M; N)= b_q(M)$, for all $q$. In particular, 
we thus recover \cite[Theorems 2.5 and 2.9]{CO}, in homological form.

\subsection{Principal ideal domain coefficients and minimal structures}
\label{subsec:pid}

Let $X$ be a minimal $CW$-complex (connected, of finite type).
Set $\pi :=\pi_1 (X)$, and denote by $C_{\bullet}(\tX)$
the $\pi$--equivariant chain complex from \S~\ref{subsec:genflag}.
By minimality,
\begin{equation}
\label{eq:ch}
C_q(\tX)= H_qX\otimes \Z \pi\, , \quad \text{for all}\quad q\, .
\end{equation}

Set $R:= \K \Z$ (the group ring of $\Z$ over a commutative field $\K$).
It is a principal ideal domain. We are going to use this fact to
prove the following $R$-analog of Corollary~\ref{cor:indarr}, 
in the context of arbitrary minimal structures.

Let $\nu : \pi \to \Z$ be an arbitrary character. 
Extend it to a change of rings, $\nu : \Z \pi \to \K \Z$.

\begin{theorem}
\label{thm:indmin}

Let  $\varphi : X \to X'$ be a homotopy equivalence, where
both $X$ and $X'$ are minimal $CW$-complexes. Use
$\varphi_{\#} :\pi_1(X)\stackrel{\sim}{\to}\pi_1(X')$
to identify fundamental groups. Then the $\K \Z$--chain complexes
$C_{\bullet}(\tX)\otimes_{\Z \pi}\K \Z$ and
$C_{\bullet}(\tX')\otimes_{\Z \pi}\K \Z$ are isomorphic,
for any change of rings homomorphism, 
$\nu : \Z \pi \to \K \Z$, as above.
\end{theorem}

The Theorem follows from the Lemma below, via 
the minimality property
\eqref{eq:ch}, and the basic homotopy invariance equation
\eqref{eq:spec}.

\begin{lemma}
\label{lem:dip}
Let $C_{\bullet}= \{ C_{q+1}\stackrel{d_{q+1}}{\longrightarrow}
C_q \}_{q\ge 0}$ and
$C'_{\bullet}= \{ C'_{q+1}\stackrel{d'_{q+1}}{\longrightarrow}
C'_q \}_{q\ge 0}$ be $R$-chain complexes, where $R$ is principal.
Assume that $C_q$ and $C'_q$ are finitely-generated free $R$-modules
of the same rank, for all $q\ge 0$, and also that
$H_q(C)$ and $H_q(C')$ are isomorphic $R$-modules, for $q\ge 0$.
Then $C_{\bullet}$ and $C'_{\bullet}$
are isomorphic $R$-chain complexes.
\end{lemma}

\begin{proof}
We may easily infer from our assumptions that the submodules
of $q$-cycles, $Z_q$ and $Z'_q$, are $R$-free, of the same rank,
$r_q$, for all $q$; likewise, the $q$-boundaries, $B_q$ and $B'_q$,
are both free of rank $s_q$, for all $q$.

Using suitable $R$-bases, the matrix of the inclusion,
$B_q \hookrightarrow Z_q$, may be put in diagonal form,
with nonzero entries, $\{ a_1, \ldots, a_{s_q} \}$, on the diagonal,
having the property that $a_1 \mid a_2\mid \cdots \mid a_{s_q}$.
Similarly, for $B'_q \hookrightarrow Z'_q$. Since 
$Z_q/B_q \simeq Z'_q/B'_q$, we infer that the elementary ideals
generated by $s \times s$ minors must be equal, for
$1\le s\le s_q$; see \cite[\S 20.2]{E}. Hence, $a_s$ and $a'_s$
differ by $R$-units, for all $s$. Therefore, 
we may find $R$-isomorphisms, 
$f_q : Z_q \stackrel{\sim}{\longrightarrow}Z'_q$, inducing
$R$-isomorphisms,
$f_q : B_q \stackrel{\sim}{\longrightarrow}B'_q$, for all
$q\ge 0$.

At the same time, we may split $d_{q+1} : C_{q+1}\surj B_q$,
for $q\ge 0$, by choosing decompositions,
$C_{q+1}= Z_{q+1}\oplus N_{q+1}$, such that
$d_{q+1}: N_{q+1}\stackrel{\sim}{\longrightarrow}B_q$.
Similarly, for $C_{\bullet}'$. Extend 
$f_q : Z_q \stackrel{\sim}{\longrightarrow}Z'_q$ to
$f_q : C_q \stackrel{\sim}{\longrightarrow}C'_q$, by
setting $f_{q\mid N_q} = (d'_q)^{-1}\circ (f_{q-1\mid B_{q-1}})
\circ d_q : N_q \stackrel{\sim}{\longrightarrow}N'_q$.
By construction, the $R$-isomorphisms $\{ f_q \}$ commute
with differentials.
\end{proof}

Theorem~\ref{thm:indmin} may be applied to non-trivial
higher homotopy groups of certain arrangements
(which are very hard to compute, in general).

Let $\A$ be an essential projective arrangement in $\PP^{r-1}$,
with complement $M$ and fundamental group $\pi$.
Assume that the cone, $\A'$, is a {\em hypersolvable}
central arrangement in $\C^r$ (see~\cite{JP1} for the definition 
and the basic properties of the hypersolvable class).
Denote by $p:= p(M)$ the order of $\pi_1$--connectivity,
introduced in ~\cite{PS}; it is a homotopy invariant of $M$,
which turns out to be combinatorial, for the hypersolvable class,
see~\cite[Corollary 4.10(1)]{PS}. Denote by $\wA$ the
{\em fiber-type deformation} of $\A$, constructed in~\cite{JP2}.

If $p=r-1$, as in \cite[Theorem 23]{DP}, then 
\cite[Theorem 18(ii)]{DP} applies to give the following
$\Z \pi$--presentation for the first nonzero higher
homotopy group of $M$:
\begin{equation}
\label{eq:pipres}
\pi_p(M(\A)) = \coker \, \{ \partial_{p+2}:
H_{p+2} M(\wA)\otimes \Z \pi \longrightarrow
H_{p+1} M(\wA)\otimes \Z \pi \}\, ,
\end{equation}
where $\partial_{p+2}$ is a boundary map from the 
$\pi$--equivariant chain complex~\eqref{eq:nottw},
associated to a Morse-theoretic minimal structure on
$M(\wA)\cong K(\pi, 1)$.

One may associate to $\A'$ a combinatorially determined
collection of positive natural numbers,
$\{ 1=d_1, d_2, \ldots, d_{\ell} \}$, called the
{\em exponents} of $\A'$, see \cite{JP1}; they
coincide with the exponents of the fiber-type central
arrangement $\wA'$ defined in \cite{FR}, see~\cite{JP2}.
One knows (\cite[Lemma 5.3]{PS}) that $\pi$ is an 
iterated semidirect product of free groups,
\begin{equation}
\label{eq:iter}
\pi = \bF_{d_{\ell}} \rtimes \cdots \rtimes \bF_{d_2}\, ,
\end{equation}
with all monodromy actions trivial on homology.

The above structural property~\eqref{eq:iter} has the following
basic practical consequence, discovered by D. Cohen and A. Suciu
in~\cite{CS2}: $K(\pi, 1)$ has a minimal structure for which
{\em all} boundary maps of the associated $\pi$--equivariant 
chain complex,
\[
\partial_q^{\text{Fox}}: H_q K(\pi, 1)\otimes \Z \pi
\longrightarrow H_{q-1} K(\pi, 1)\otimes \Z \pi \, ,
\]
may be {\em explicitly} computed, by Fox differential calculus.

\begin{corollary}
\label{cor:piexpl}

Let $\A$ be an essential arrangement in $\PP^{r-1}$, with
hypersolvable cone, $\A'$. Set $M:= M(\A)$, $\pi:= \pi_1(M)$,
$p:= p(M)$. Assume that $p=r-1$. Let $\nu : \Z \pi \to \K \Z$
be the change of rings associated to an arbitrary character,
$\nu : \pi \to \Z$ (where $\K$ is a commutative field). Then :
\begin{equation}
\label{eq:fox}
\pi_p(M(\A))\otimes_{\Z \pi}\K \Z = \coker \,
\{ \partial_{p+2}^{\text{Fox}} \otimes_{\Z \pi}\K \Z :
(\K \Z)^{b_{p+2}(\pi)}\longrightarrow
(\K \Z)^{b_{p+1}(\pi)} 
\}\, ,
\end{equation}
as $\K \Z$--modules, where the Betti numbers of $\pi$ are
determined by the exponents of $\A'$, and $\partial_{\bullet}^{\text{Fox}}$
is explicitly computed from~\eqref{eq:iter}. 
\end{corollary}

\begin{proof}
The Poincar\' e polynomial of $K(\pi, 1)\cong M(\wA)$ is
$P(T)= \prod_{i=2}^{\ell} (1+ d_i T)$; see~\cite{FR}.
From \eqref{eq:pipres}, we infer that $\pi_p (M)\otimes_{\Z \pi}
\K \Z = \coker \, \{ \partial_{p+2}\otimes_{\Z \pi}
\K \Z \}$. By Theorem~\ref{thm:indmin}, we may replace
$\partial_{p+2}$ by $\partial_{p+2}^{\text{Fox}}$, to arrive 
at \eqref{eq:fox}, as asserted.
\end{proof}

\section{A Milnor fiber obstruction to relative minimality}
\label{sec:relmin}

\subsection{Absolute minimality and $H_1$--bases}
\label{subsec:abs}

We come back to Definition~\ref{def:minad} 
from the preceding section.
We will add more information, related to distinguished $H_1$--bases.
This will be needed for certain twisted homology
computations (such as those related to the homology of
Milnor fibers).

Let $\{ H_0, H_1,\ldots ,H_n \}$ be the hyperplanes of $\A$,
in $\PP (U)$. The {\em meridians} associated to the hyperplanes,
$\{ \mu_i \in \pi_1(M) \}_{0\le i\le n}$, give a collection of 
well-defined elements of $H_1 M$. Considering $H_0$ as a 
distinguished hyperplane, we obtain in this way a
distinguished $\Z$--basis of $H_1 M$,
\begin{equation}
\label{eq:mbase}
\{ [\mu_i] \}_{1\le i\le n} \, .
\end{equation}
Let now $\{ c_i \}_{1\le i\le n}$ be the $1$--cells of a
minimal $CW$-structure on $M$, $X$. They provide a
distinguished $\Z$--basis of $H_1 X$, 
\begin{equation}
\label{eq:cbase}
\{ [c_i] \}_{1\le i\le n}\, .
\end{equation}

\begin{definition}
\label{def:mark}

Let $X$ be a minimal complex. A homotopy equivalence,
$\varphi : X\to M$, {\em respects $H_1$--markings},
if it takes the basis~\eqref{eq:cbase} to the basis
\eqref{eq:mbase}.
\end{definition}

Every arrangement complement, $M$, has such a {\em marked}
minimal structure. Indeed, take an $\A$--generic flag $\F$,
as in Definition~\ref{def:flag}, and consider 
an $\F$--compatible minimal structure, as in 
Definition~\ref{def:minad}. It is constructed inductively, see
\cite[Section 4]{DP}. At the first nontrivial step,
$X^{(1)}$ is a wedge of $n$ circles, and $F_1 M$ is
$\PP1 \setminus \{ n+1 \,\, \text{points}\,  \}$.
Obviously, we may start with a homotopy equivalence,
$\varphi_1$, which preserves the canonical $H_1$--bases,
and then proceed by induction.

\subsection{The relative minimality problem}
\label{subsec:rel}

Assume now that $\A$ is an $\LL_k(\wA)$--generic section
of $\wA$ (in the sense from~\cite[(1)]{DP}), with
$k\ge 1$, where $\wA$ is an essential arrangement of
$n+1$ hyperplanes in $\PP (V)$. Our basic idea in~\cite{DP}
was to extract from $\wA$ homotopy information on $M(\A)$.
The key tool is provided by \cite[Proposition 14]{DP},
which says that one may replace, up to homotopy,
the inclusion, $j : M(\A)\hookrightarrow M(\wA)$, by
a cellular map between minimal complexes, $f : X\to Y$,
with the property that
\begin{equation}
\label{eq:ksk}
f_{\mid X^{(k)}} = \id \,\, .
\end{equation}

The {\em relative minimality} problem we have in mind is
the following. Start with a marked minimal structure on
$M(\wA)$, $Y$, as in Definition~\ref{def:mark}. Let $\A$
be an essential $\LL_k$--generic section of $\wA$, with
$k\ge 1$. Can one replace $j$ by $f$, as in~\eqref{eq:ksk} above, 
in such a way that moreover
\begin{equation}
\label{eq:sub}
f_{\mid X^{(k+1)}}\quad =\quad \text{subcomplex inclusion}\quad ?
\end{equation}

When $\A$ has hypersolvable cone and $p:= p(M(\A))=
\rank \,(\A') -1$, as in Corollary~\ref{cor:piexpl},
this can be done: actually one may take $X =Y^{(p)}$,
and $f=$ inclusion; see~\cite[Theorems 18 and 23]{DP}.
In general, the (easily checked) fact that $H_*j$ is
a split injection, together with heuristic morsification
arguments, seem to indicate that the answer to question
\eqref{eq:sub} should be yes.

Surprisingly enough, the answer turns out to be no, in
general, and the homology of the Milnor fiber of $\A$
comes into play, at this point.

\subsection{Twisted homology and Betti numbers of Milnor fibers}
\label{subsec:bigcof}

Let $q:= \prod_{i=0}^n \alpha_i$ be a defining equation of 
the central arrangement in $U$ associated to $\A$. Let
$F:= q^{-1}(1)$ be the Milnor fiber of $\A$. We are going to
recall from~\cite[Corollary 1.5]{CS1} the twisted homology
decomposition of $H_*(F; \C)$.

Set $u:= \exp \,(\frac{2\pi \sqrt{-1}}{n+1})$. For
$0\le t\le n$, denote by $L_t$ the rank one $\C$--local
system on $M:= M(\A)$ (alias, the abelian representation,
$L_t : H_1 M \to \C^*$, of $\pi_1(M)$) given by
\begin{equation}
\label{eq:eqm}
L_t([\mu_i]) = u^t \, ,\quad \text{for} \quad 1\le i \le n
\end{equation}
(where $\{ [\mu_i] \}$ is the $\Z$--basis~\eqref{eq:mbase}). Set
\begin{equation}
\label{eq:defspl}
b_s^t \, (F):= \dim_{\C} H_s(M; L_t)\, .
\end{equation}
Then:
\begin{equation}
\label{eq:splitf}
b_s (F)= \sum_{t=0}^n b_s^t \, (F)\, , \quad \text{for all}\quad s\, .
\end{equation}

\subsection{The Milnor fiber obstruction}
\label{subsec:obstr}

Let $\A$ be an arbitrary essential arrangement in $\PP (U)$,
with defining equation $q$ and Milnor fiber $F$, as in 
\S~\ref{subsec:bigcof}. The linear forms $\{ \alpha_i \}_{0\le i\le n}$
define a linear embedding
\begin{equation}
\label{eq:embess}
j : U\hookrightarrow V:= \C^{n+1}\, ,
\end{equation}
which enables us to view $\A$ as an $\LL_1$--generic section of
the Boolean arrangement $\wA$ (with defining equation
$\widehat{q}:= \prod_{i=0}^n z_i$).

Set $Y:= (S1)^{\times n}$, endowed with the canonical minimal
structure of the $n$-torus. Plainly, there is a marked homotopy
equivalence, $\varphi_Y : Y\to M(\wA)$. Assume that
$j: M(\A)\hookrightarrow M(\wA)$ has the homotopy type of
a cellular map between minimal complexes, $f: X\to Y$, 
with the property that $f_{\mid X^{(1)}}= \id $, as in \eqref{eq:ksk},
and $f_{\mid X^{(2)}}=$ inclusion, as in~\eqref{eq:sub}. (Note that
$\varphi_X : X\to M(\A)$ will also be a marked homotopy equivalence,
since, obviously, $j$ and $f$ respect $H_1$--markings.)   

\begin{theorem}
\label{thm:best}

Let $\A$ be an essential arrangement of $n+1$ hyperplanes,
with Milnor fiber $F$. If the relative minimality problem
\eqref{eq:sub}, where $k=1$, has a positive answer, for $\A$
and the Boolean arrangement $\wA$, then $n$ divides $b_1(F)$.
\end{theorem}

\begin{proof}
Since $b_10(F)= b_1(M(\A))= n$, it will be enough to
show that $b_1^t(F)$ is independent of $t$, for $1\le t \le n$;
see~\eqref{eq:splitf}. As $\varphi_X$ is a marked homotopy equivalence,
the twisted Betti numbers from~\eqref{eq:defspl} may be computed on
$X$, using the cellular $H_1$--basis \eqref{eq:cbase} in \eqref{eq:eqm}, 
via the basic specialization formula~\eqref{eq:spec}.

Our hypothesis~\eqref{eq:ksk} on $f_{\mid X^{(1)}}$ readily implies that
$f_{\#}: \pi := \pi_1(X) \to \pi_1(Y)= \Z^n$ is the
abelianization map. The cellular map $f: X\to Y$ lifts to a
$\Z f_{\#}$--linear chain map, 
\[
\{ \widetilde{f_s} : H_s X\otimes \Z \pi\longrightarrow
H_s Y\otimes \Z \Z^n \}_{s\ge 0}\, ,
\]
between the equivariant chain complexes of the universal covers,
\[
C_{\bullet}(\tX):= \{ H_{s+1} X\otimes \Z \pi 
\stackrel{d_{s+1}}{\longrightarrow}
 H_{s} X\otimes \Z \pi \}_{s\ge 0}\, ,
\]
and
\[
C_{\bullet}(\tY):= \{ H_{s+1} Y\otimes \Z \Z^n 
\stackrel{\partial_{s+1}}{\longrightarrow}
 H_{s} Y\otimes \Z \Z^n \}_{s\ge 0}\, .
\]

Our main hypothesis~\eqref{eq:sub} implies that
\begin{equation}
\label{eq:tildef}
\widetilde{f_s} = H_s f \otimes \Z f_{\#} \, ,
\quad \text{for} \quad s\le 2\, ,
\end{equation}
with $H_{\le 2} f$ monic.

For $s\le 1$, we may thus tensor the commutative squares
\begin{equation*}
\xymatrix  @C5pc {
H_{s+1} X\otimes\Z \pi  \ar[r]^{d_{s+1}}
\ar[d]^{H_{s+1} f\otimes \Z f_{\#}}
& H_s X\otimes\Z \pi  \ar[d]^ {H_s f\otimes \Z f_{\#}}
\\
H_{s+1} Y\otimes\Z \Z^n  \ar[r]^{\partial_{s+1}}
& H_s Y\otimes\Z \Z^n
}
\end{equation*}
(see~\eqref{eq:tildef}) with $\Z \Z^n$ over $\Z \pi$, via
$\Z f_{\#}$, to get commuting squares
\begin{equation*}
\xymatrix  @C5pc {
H_{s+1} X\otimes\Z \Z^n  \ar[r]^{d_{s+1}^{\text{ab}}}
\ar[d]^{H_{s+1} f\otimes \id}
& H_s X\otimes\Z \Z^n  \ar[d]^ {H_s f\otimes \id}
\\
H_{s+1} Y\otimes\Z \Z^n  \ar[r]^{\partial_{s+1}}
& H_s Y\otimes\Z \Z^n
}
\end{equation*}

For any $1\le t\le n$, we may further specialize to $\C$, via the
representation~\eqref{eq:eqm}. In this way, we get commutative squares,
\begin{equation}
\label{eq:comm}
\xymatrix  @C5pc {
H_{s+1} X\otimes\C  \ar[r]^{d_{s+1}^{\text{ab}}(u^t)}
\ar[d]^{H_{s+1} f\otimes \id}
& H_s X\otimes\C  \ar[d]^ {H_s f\otimes \id}
\\
H_{s+1} Y\otimes\C  \ar[r]^{\partial_{s+1}(u^t)}
& H_s Y\otimes\C
}
\end{equation}
(for $s\le 1$), where the vertical maps are injective, 
and independent of $t$. Note that the upper chain complex 
from~\eqref{eq:comm} computes $b_1^t(F)$; 
see~\eqref{eq:defspl} and~\eqref{eq:spec}.

At the same time, eye-inspection of the well-known explicit
formula for $\{ \partial_{s+1} \}_s$ (see e.g. 
\cite[(10)]{DP}) reveals that 
$\partial_{s+1}(u^t)= (u^{-t} -1)\cdot \partial'_{s+1}$,
where the differential $\partial'_{s+1}$ is independent of $t$.
By~\eqref{eq:comm} above, $b_1^t(F)$ is therefore 
independent of $t$, for
$0<t\le n$.
\end{proof} 

\begin{example}
\label{ex:no}

In Example 5.1 from~\cite{CS1}, $n=5$, $b_1^t=0$, for
$t=1,3,5$, and $b_1^t=1$, for $t=2,4$; hence, $b_1(F)=7$.
In Example 5.4 from~\cite{CS1}, $n=8$, $b_1^t=0$, for
$t=1,2,4,5,7,8$, and $b_1^t=1$, for $t=3,6$; hence, $b_1(F)=10$.
By Theorem~\ref{thm:best}, the relative minimality problem
\eqref{eq:sub} has a negative answer, in both cases.
\end{example}

\section{Twisted homology with resonant coefficients}
\label{sec:tord}

Let $\A$ be an essential, proper, $\LL_k(\wA)$--generic
section, with $k\ge 2$, of an essential 
{\em aspherical} arrangement, $\wA$. 
(When we say `proper', we want to exclude the trivial case,
$\wA =\A$.) Set $M:= M(\A)$, $\pi :=\pi_1(M)$, and
$p:= p(M)$.

We know that $j: M(\A)\hookrightarrow M(\wA)$ has
the homotopy type of a cellular map between minimal 
complexes, $f: X\to Y$,
with the property that
\begin{equation}
\label{eq:htpyf}
f_{\mid X^{(p)}} =\id \, ,
\end{equation}
where $2\le k\le p<\infty$. In particular, $Y$ is a
$K(\pi, 1)$. See the discussion preceding Theorem 16~\cite{DP}.

Our goal in this section is to use~\eqref{eq:htpyf} above,
to perform various twisted homology computations on $M$,
in terms of $\pi_1(M)$. We will {\em not}
impose any kind of `nonresonance' conditions on the coefficients.

\subsection{Computations in the $\LL$--generic range}
\label{subsec:lgen}

Here the coefficients will be quite general, but our method will
give results only in the {\em $\LL$--generic range}, that is,
up to $H_{p-1}(M; L)$.

\begin{theorem}
\label{thm:range}

Let $\A$ be an essential, proper, $\LL_k(\wA)$--generic
section, with $k\ge 2$, of an essential 
aspherical arrangement, $\wA$. 
Set $M:= M(\A)$, $\pi :=\pi_1(M)$, and
$p:= p(M)$. 
\begin{enumerate}
\item \label{a1}
Let $L$ be an arbitrary local system on $M$. Then:
\[
H_q(M; L)= \tor_q^{\Z \pi} (\Z, L)\, ,\quad \text{for} \quad q<p \, .
\]
\item \label{a2}
If $\wA$ is {\em fiber-type}, 
then the $\tor^{\Z \pi}$--groups from Part~\eqref{a1} may
be explicitly computed using the Fox $\Z \pi$--resolution 
of $\Z$ from~\cite{CS2}.
\end{enumerate}
\end{theorem}

\begin{proof}
Part~\eqref{a1}. By~\eqref{eq:spec} and basic 
homotopy formula~\eqref{eq:htpyf},
\[ 
H_{<p} (M; L)= H_{<p} (C_{\bullet}(\tY)\otimes_{\Z \pi} L)\, .
\]
The identification of $ H_{<p} (C_{\bullet}(\tY)\otimes_{\Z \pi} L)$
with $\tor_{<p}^{\Z \pi} (\Z, L)$ comes now from the fact that
the $\pi$--equivariant chain complex $C_{\bullet}(\tY)$ is a free
{\em $\Z \pi$--resolution} of $\Z$, since $Y$ is a $K(\pi, 1)$.

Part~\eqref{a2}. If $\wA$ is fiber-type, then 
$\pi =\pi_1(M(\wA))$ is an iterated semidirect product of
free groups, with trivial monodromy actions on homology,
as in \eqref{eq:iter}. Therefore, the Fox calculus free
$\Z \pi$--resolution of $\Z$ from~\cite{CS2} may be used as well
to compute the  $\tor^{\Z \pi}$--groups from Part~\eqref{a1}.
\end{proof}

As an illustration of Theorem~\ref{thm:range}, we may offer
the following simple, very explicit, class of examples. Let
$\A$ be an essential projective arrangement of $n+1$ hyperplanes,
with associated central arrangement, $\A'$. Given a subarrangement,
$\B' \subset \A'$, denote by $\mid \B' \mid$ the number 
of hyperplanes of $\B'$.

Define $c(\A):= \infty$, if $\A'$ is independent. Otherwise, set
\begin{equation}
\label{eq:defc}
c(\A):= \min \,\, \{ \mid \B' \mid \, \mid \, 
\B' \subset \A' \quad \text{is dependent}\, \}\, .
\end{equation}
Obviously, $c(\A)\ge 3$.

\begin{corollary}
\label{cor:ccomb}

Let $\A$ be an essential arrangement of $n+1$
hyperplanes in $\PP (U)$. Set $M:= M(\A)$,
$\pi:= \pi_1(M)$, $p:= p(M)$, and $c:= c(\A)$.
\begin{enumerate}
\item \label{b2}
If $c<\infty$, then $\A$ is a proper $\LL_{c-2}$--generic
section of an essential Boolean arrangement of $n+1$
hyperplanes.
\item \label{b3}
If $c>3$, then $\pi= \Z^n$ and $p= c-2$.
\item \label{b4}
Assume $3<c<\infty$. Let $L$ be a left $\Z \Z^n$--module.
Define a chain complex,
\begin{equation}
\label{eq:cl}
C_{\bullet}(\Z^n; L):= \{ \bigwedge^s (x_1, \ldots,x_n)
\otimes_{\Z} L \stackrel{\partial_s}{\longrightarrow}
\bigwedge^{s-1}(x_1, \ldots,x_n)\otimes_{\Z} L \}_{s\ge 1}\, ,
\end{equation}
by setting
\begin{equation}
\label{eq:defcl}
\partial_s (x_{i_1}\wedge \cdots \wedge x_{i_s}\otimes v) 
:= \sum_{r=1}^s (-1)^{r-1}
x_{i_1}\wedge \cdots \widehat{x_{i_r}}\cdots \wedge x_{i_s}
\otimes (x_{i_r}^{-1} -1) v\, ,
\end{equation}
where $\{ x_1,\ldots, x_n \}$ denotes the standard basis of $\Z^n$. Then:
\[
H_q(M; L)= H_q(C_{\bullet}(\Z^n; L))\, ,\quad \text{for} \quad
q<c-2 \, .
\]
In particular, the $L$--twisted homology of $M(\A)$ in the
$\LL$--generic range depends only on $L$, and the 
combinatorics of $\A$. 
\end{enumerate}
\end{corollary}

\begin{proof}
Part~\eqref{b2}. Equation~\eqref{eq:embess} from
\S~\ref{subsec:obstr} shows that $\A$ is an $\LL_1$--generic section
(proper, since $c<\infty$) of the required Boolean arrangement, $\wA$.
This section may be easily seen to be actually $\LL_{c-2}$--generic,
by resorting to the definitions: see~\cite[(1)]{DP} for 
$\LL_k$--genericity, and~\eqref{eq:defc} above for $c(\A)$.

Part~\eqref{b3}. If $c=\infty$, then plainly $\A$ itself is Boolean,
$\pi =\Z^n$ and $p=\infty$. Assume then that
$3<c<\infty$. By~\cite[Proposition 14]{DP}, $\pi =\Z^n$ and
$M(\wA)$ is a $K(\pi, 1)$. Moreover,
$k(\wA, U)=p(\wA, U)$, according to~\cite[Proposition 15]{DP}.
One may check the equality $k(\wA, U)=c-2$, directly from the
definitions (\cite[(3)]{DP} and~\eqref{eq:defc} respectively).
Likewise, the equality $p(\wA, U)= p(M(\A))$ follows 
from the definitions
(\cite[(4)]{DP} and~\cite[p.73]{PS} respectively). Finally,
$p=c-2$, as asserted.

Part~\eqref{b4}. Parts~\eqref{b2}--\eqref{b3} enable us to compute
$H_{<c-2}(M; L)$ as in Theorem~\ref{thm:range}~\eqref{a1}. 
The equality $H_q(C_{\bullet}(\Z^n; L))= \tor_q^{\Z \Z^n}(\Z, L)$
follows at once, by using the standard $\Z \Z^n$--resolution of $\Z$
\cite[(10)]{DP} to compute $\tor^{\Z \Z^n}$--groups.
\end{proof}

\subsection{A fiber-type general position framework for complete
computations}
\label{subsec:genhat}

Now we are going to focus on the {\em fiber-type general position} class
from Corollary~\ref{cor:piexpl}. By definition, an essential arrangement
$\A$ belongs to this class if it has hypersolvable cone, $\A'$, and
$p:= p(M(\A))=r-1$, where $r:= \rank \,(\A')$. We know that $p\ge 2$, and
that this class coincides with the proper $\LL_{r-1}$--generic sections,
$r\ge 3$, of essential fiber-type arrangements, $\wA$; see
the discussion preceding Theorem 23~\cite{DP}.

When $\wA$ is Boolean, the above definition corresponds to
the {\em general position} arrangements, intensively studied
since Hattori's pioneering work~\cite{H}; see~\cite[Remark 19]{DP}.
This explains both our teminology, and our interest in 
the fiber-type general position class.

Let $\K$ be a commutative field. We will treat
two types of local coefficients on $M:= M(\A)$.
The first type consists of $PID$
coefficients, $R:= \K \Z$, coming from characters,
$\nu : \pi\to \Z$, as in \S~\ref{subsec:pid}. The second type
consists of $\pi$--modules $N:= \K^d$, coming from
finite-dimensional $\K$--representations, 
$\rho : \pi\to GL(d; \K)$. In both cases, our choice was guided by
the desire to have an Euler characteristic argument at hand,
to obtain complete computations.

\begin{theorem}
\label{thm:compl}

Let $\A$ be an essential arrangement in $\PP^{r-1}$, with
$r\ge 3$. Set $M:= M(\A)$ and $\pi:= \pi_1(M)$. Assume that $\A$ is 
a proper, $\LL_{r-1}(\wA)$--generic section, of an essential aspherical
arrangement, $\wA$. Let $\K$ be a commutative field.
\begin{enumerate} 
\item \label{c1}
Set $R:= \K \Z$, endowed with the left $\Z \pi$--module structure
coming from a character, $\nu : \pi\to \Z$. Then:
\[
H_q(M(\A); R)=
\left\{
\begin{array}{cll}
\tor_q^{\Z \pi} (\Z, R) \, , & \text{for} & q< r-1\, ;\\
\text{$R$-free, of rank}\, =(-1)^{r-1}[\chi(M(\A))- \kappa_{r-1}]\, ,
& \text{for} & q=r-1\, ;\\
0\, , & \text{for} & q>r-1\, , 
\end{array}
\right.
\]
where $\kappa_{r-1}= \sum_{q=0}^{r-2} (-1)^q \rank_R \tor_q^{\Z \pi}
(\Z, R)$.
\item \label{c2}
Set $N:= \K^d$,  endowed with the left $\Z \pi$--module structure
coming from a representation, 
$\rho : \pi\to GL(d; \K)$. Then:
\[
H_q(M(\A); N)=
\left\{
\begin{array}{cll}
\tor_q^{\Z \pi} (\Z, N) \, , & \text{for} & q< r-1\, ;\\
\text{of $\K$-dim}\, =(-1)^{r-1}[d\cdot \chi(M(\A))- \kappa_{r-1}]\, ,
& \text{for} & q=r-1\, ;\\
0\, , & \text{for} & q>r-1\, , 
\end{array}
\right.
\]
where $\kappa_{r-1}= \sum_{q=0}^{r-2} (-1)^q \dim_{\K} \tor_q^{\Z \pi}
(\Z, N)$.
\item \label{c3}
If $\wA$ is fiber-type (e.g., if $\A$ has hypersolvable cone,
and $p(M(\A))=r-1$), then all $\tor^{\Z \pi}$--computations from Parts~
\eqref{c1} and~\eqref{c2} may be done explicitly, with the aid of
the Fox $\Z \pi$--resolution of $\Z$ from~\cite{CS2}.
\end{enumerate}
\end{theorem}

\begin{proof}
Set $p:= p(M)$. We know that $p=r-1$, by~\cite[Theorem 18(i)]{DP}.
For $q<r-1$, our assertions from Parts~\eqref{c1}--\eqref{c2}
follow then directly from Theorem~\ref{thm:range}\eqref{a1}.
If $q>r-1$, then plainly $H_q(M; L)=0$ ($L$ arbitrary), for
dimensional reasons. It remains to settle the case $q=r-1$,
for $R$ and for $N$.

Denote by $X$ (respectively $Y$) the minimal $CW$-structure on
$M(\A)$ (respectively on $M(\wA)$). We infer from~\eqref{eq:htpyf} 
that $X= Y^{(r-1)}$, where $Y$ is a $K(\pi, 1)$. Therefore
(see~\eqref{eq:spec}), 
$H_{r-1}(M; L)= \ker \, \{ \partial_{r-1}\otimes_{\Z \pi}
L : H_{r-1} Y\otimes_{\Z} L \longrightarrow 
H_{r-2} Y\otimes_{\Z} L \}$, for arbitrary $L$, where
$\{ \partial_{r-1} \}_{r\ge 2}$ 
denotes the differential 
of the $\pi$--equivariant chain complex $C_{\bullet}(\tY)$.

In Part~\eqref{c1}, we may use an Euler characteristic argument for
the finite, $R$--free, $R$--chain complex
$C_{\bullet}(\tY^{(r-1)})\otimes_{\Z \pi} R$,
to infer that $H_{r-1}(M; R)$ is $R$--free and
\[
\chi (M)= (-1)^{r-1} \rank_R H_{r-1}(M; R) + \kappa_{r-1} \, ,
\]
as asserted.

In Part~\eqref{c2}, we may apply the same argument to the finite
$\K$--chain complex $C_{\bullet}(\tY^{(r-1)})\otimes_{\Z \pi} N$,
to get
\[
d\cdot \chi (M)= (-1)^{r-1} \dim_{\K} H_{r-1}(M; N) + \kappa_{r-1} \, ,
\]
which verifies our claimed formula.

As for Part~\eqref{c3}, we may use the same argument as 
in the proof of Theorem~\ref{thm:range}\eqref{a2}.
\end{proof}

\section{Character--abelianizations of higher homotopy groups}
\label{sec:ab}

\subsection{Nonresonant abelianization and combinatorics}
\label{subsec:nonres}

We continue to study the homotopy properties of the complement,
for arrangements belonging to the fiber-type general position class
from Corollary~\ref{cor:piexpl}. Within this class,
we have explained, in \S~\ref{subsec:pid}, how to construct 
an explicit presentation matrix for
$\pi_p (M)\otimes_{\Z \pi} \K \Z$, for an arbitrary 
character, 
$\nu : \pi\to \Z$.

Apriori, only the size of the presentation matrix
$\partial_{p+2}^{\text{Fox}}\otimes_{\Z \pi} \K \Z$
from Corollary~\ref{cor:piexpl} is {\em combinatorially determined}.
In this subsection, we aim at enlarging the dictionary
`topology $\leftrightarrow$ combinatorics', along 
the lines from~\cite{DP}; see especially Theorem 23 therefrom.

More precisely, we would like to identify more numerical 
invariants of $\pi_p(M(\A))$ which are determined
by the combinatorics of $\A$. A natural candidate is provided by
$\rank_R \pi_p(M)\otimes_{\nu} R$, where $R:= \K \Z$. 
To derive a combinatorial formula for 
the aforementioned $R$--rank,
we are led to impose a certain `nonresonance'
condition on the local system $\nu$.

This condition in turn is inspired from 
the powerful vanishing result,
for $\C$--local systems of finite rank, proved in~\cite{CDO}.
To make the appropriate formal definition, we recall that
there is a combinatorially defined subset of edges,
$\D (\A)\subset \LL (\A)$, called {\em dense} edges,
for any arrangement $\A$; see~\cite{STV}.

\begin{definition}
\label{def:nr}

Let $\A = \{ H_0, H_1,\ldots, H_n \}$ be an essential arrangement,
with fundamental group $\pi$. A character,
$\nu : \pi\to \Z$ (that is, a collection 
$\{ \gamma_i:= \nu (\mu_i)\in \Z \}_{0\le i\le n}$,
such that $\sum_{i=0}^n \gamma_i =0$)
is called {\em $\A$--nonresonant} if
\begin{equation}
\label{eq:anr}
\sum_{H_i \supset S} \gamma_i \ne 0 \, , \,
\forall S\in \D(\A) \quad \text{such that} \quad 
S\subset H_0 \, .
\end{equation}
\end{definition}

Given $\A$, it is straightforward to check that
$\A$--nonresonant characters always exist.

\begin{theorem}
\label{thm:rcomb}

Let $\A$ and $\wA$ be essential arrangements, in 
$\PP^{r-1}$ and $\PP^{m-1}$ respectively. Assume that
$\wA$ is aspherical, and $\A$ is a proper,
$\LL_{r-1}(\wA)$--generic section of $\wA$, 
with $r\ge 3$. Set $M:= M(\A)$ and $\pi:= \pi_1(M)$.
\begin{enumerate}
\item \label{d1}
$\pi_{r-1}(M(\A))$ is the first higher non-vanishing
homotopy group of $M$.
\item \label{d2}
Let $\nu : \pi\to \Z$
be $\wA$--nonresonant (in the sense of Definition~\ref{def:nr} above).
Set $R:= \C \Z$. Then:
\[
\rank_R \pi_{r-1}(M(\A))\otimes_{\nu} R=
\left\{
\begin{array}{cll}
(-1)^{r-1} \chi (M(\A))\, , & \text{if} & r+1<m \, ;\\
b_r (\pi)\, , & \text{if} & r+1 =m \, .
\end{array}
\right.
\]
\item \label{d3}
If $\A$ has hypersolvable cone, $p(M(\A))=r-1$, and
$\wA$ is the fiber-type deformation of $\A$ from~\cite{JP2},
then $b_r(\pi)$ from the second case in Part~\eqref{d2} equals
$\prod_{j=1}^{\ell} d_j$, where $\{ d_j \}$ are 
the exponents of the cone of $\A$. In particular,
$\rank_R \pi_{r-1}(M(\A))\otimes_{\nu} R$ is
combinatorially determined, if $\nu$ is $\wA$--nonresonant.
\end{enumerate}
\end{theorem}

\begin{proof}
Part~\eqref{d1} follows directly from~\cite[Theorem 18]{DP}.
Let $Y$ be the minimal $CW$-structure on $M(\wA)$
from~\cite[Theorem 18(ii)]{DP}, with associated 
$\pi$--equivariant chain complex,
$C_{\bullet}(\tY):= \{ \partial_q : H_q Y\otimes \Z \pi
\longrightarrow H_{q-1} Y\otimes \Z \pi \}_{q\ge 1}$. We
infer, from the $\Z \pi$--resolution~\cite[(11)]{DP} 
of $\pi_{r-1}(M)$, that
\begin{equation}
\label{eq:rpres}
\pi_{r-1}(M)\otimes_{\nu} R =\coker \, 
\{ H_{r+1} Y\otimes_{\Z} R 
\stackrel{\partial_{r+1}\otimes_{\nu} R}{\longrightarrow}
H_r Y\otimes_{\Z} R \} \, ,
\end{equation}
for arbitrary $\nu$.

Part~\eqref{d2}. If $r+1=m$, then
$\pi_{r-1}(M)\otimes_{\nu} R =H_r Y\otimes R$
is $R$--free, with $R$--rank equal to $b_r(\pi)$, 
as asserted, for arbitrary $\nu$
(since $Y$ is a $K(\pi, 1)$). Assume then that $r+1<m$.

In this case, $\wA$--nonresonance comes into play, in the
following way. As a preliminary remark, note that the chain complex,
$C_{\bullet}(\tY)\otimes_{\Z \pi} \C_t :=
\{ \partial_{q+1}\otimes_{\nu} R (t) :
H_{q+1} Y\otimes \C \longrightarrow H_q Y \otimes \C \}_{q\ge 0}$,
obtained by further specializing $1\in \Z$ to $t\in \C^*$,
computes the homology of $M(\wA)$ with coefficients in the
appropriately defined rank one $\C$-local system, $\C_t$;
see~\eqref{eq:spec}.

Therefore, Theorem 1 and Lemma 2 from~\cite{CDO} together
imply that
\begin{equation}
\label{eq:vanish}
H_{\ne m-1} (C_{\bullet}(\tY)\otimes_{\Z \pi} \C_t)=0 \, ,
\end{equation}
as soon as
\begin{equation}
\label{eq:tgen}
t^{\sum_{H_i\supset S}\gamma_i}\ne 1\, ,
\forall S \in \D (\wA)\quad \text{such that} \quad
S \subset H_0 \, .
\end{equation}
Comparing~\eqref{eq:tgen} with~\eqref{eq:anr},
we infer from $\wA$--nonresonance that the vanishing property
\eqref{eq:vanish} holds, for generic $t$.

Recall now from~\eqref{eq:htpyf} that $Y^{(r-1)}$ is
a minimal $CW$-structure for $M(\A)$. An Euler
characteristic argument,
applied to the $\C$--chain complex
$C_{\bullet}(\tY^{(r-1)})\otimes_{\Z \pi} \C_t$, 
for generic $t$, provides the equality
\begin{equation}
\label{eq:chi}
(-1)^{r-1} \chi(M(\A))= \dim_{\C}\ker \,
\{ \partial_{r-1}\otimes_{\nu} R (t) \} \, .
\end{equation}

Since $r<m-1$, \eqref{eq:vanish} also implies that
\begin{equation}
\label{eq:kc}
\ker \, \{ \partial_{r-1}\otimes_{\nu} R (t) \}=
\coker \, \{ \partial_{r+1}\otimes_{\nu} R (t) \}\, ,
\end{equation}
for generic $t$. Using equations~\eqref{eq:chi},
\eqref{eq:kc} and~\eqref{eq:rpres}, we arrive at
\begin{equation}
\label{eq:gata}
\dim_{\C}(\pi_{r-1}(M)\otimes_{\nu} R)\otimes_R \C_t =
(-1)^{r-1} \chi(M(\A))\, ,
\end{equation}
for generic $t$.

At the same time, it is well-known that 
$\rank_R N =\dim_{\C} (N\otimes_R \C_t)$, generically,
for a finite $R$--module $N$. With this remark,
equation~\eqref{eq:gata} finishes the proof of Part~\eqref{d2}.

Part~\eqref{d3}. In the hypersolvable case, it is enough to
recall from~\cite{FR} that the Poincar\' e polynomial
of $M(\wA)\cong K(\pi, 1)$ is $\prod_{j=2}^{\ell} (1+ d_j T)$,
where $\ell =m$, and $d_1=1$; see also our discussion preceding
Corollary~\ref{cor:piexpl} from \S\ref{subsec:pid}.
\end{proof}

\begin{remark}
\label{rem:tors}

Let $\A$ be an arbitrary essential arrangement in $\PP^{r-1}$,
with fundamental group $\pi$. Set $R:= \C \Z$. 
Let $\nu : \pi\to \Z$ be
$\A$--nonresonant. Then it is not difficult to use the
same kind of arguments as in the proof of Theorem 4.2
from~\cite{DN}, to obtain that
\[
H_q(M(\A); R)=
\left\{
\begin{array}{cll}
\text{$R$--torsion} \, , & \text{for} & q< r-1\, ;\\
\text{$R$-free, of rank}\, =(-1)^{r-1} \chi(M(\A))\, ,
& \text{for} & q=r-1\, ;\\
0\, , & \text{for} & q>r-1\, . 
\end{array}
\right.
\] 
\end{remark}

\subsection{General abelianization}
\label{subsec:genab}

We may drop the nonresonance restriction from Theorem
\ref{thm:rcomb}\eqref{d2}. We obtain the following
(apriori, non-combinatorial) formula for the rank
of an arbitrary character--abelianization, in terms of the
fundamental group and the Euler characteristic.

\begin{theorem}
\label{thm:rarb}

Let $\A$  be an essential arrangement in 
$\PP^{r-1}$, where $r\ge 3$. Assume that $\A$ is a proper,
$\LL_{r-1}(\wA)$--generic section of
an essential aspherical arrangement, $\wA$. 
Set $M:= M(\A)$ and $\pi:= \pi_1(M)$.
Let $\nu : \pi\to \Z$ be an arbitrary character.
Denote by $R:= \K \Z$ the associated local system on $M$,
where $\K$ is any commutative field. Then:
\[
\rank_R \pi_{r-1}(M(\A))\otimes_{\nu} R=
(-1)^{r-1}[ \chi (M(\A))
- \sum_{q=0}^r (-1)^q \rank_R \tor_q^{\Z \pi} (\Z, R)] \, .
\]
\end{theorem} 

\begin{proof}
The homotopy formula~\eqref{eq:htpyf} says that
$Y^{(r-1)}$ is a minimal $CW$-structure on $M(\A)$, where
$Y$ is a minimal structure on $M(\wA)$.

We are going to use the K\" unneth spectral sequence from
\cite[Theorem XII.12.1]{ML}, arising from the free right
$\Z \pi$--chain complex, $C_{\bullet}(\tY^{(r-1)})$, and
the change of rings map, $\nu : \Z \pi\to R$.

The $E2$--term, $E_{st}^2 = \tor_s^{\Z \pi}(H_t \tY^{(r-1)}, R)$,
is possibly nonzero only for $t=0$ (where $H_t \tY^{(r-1)}= \Z$),
and for $t=r-1$ (where $H_t \tY^{(r-1)}=\pi_{r-1}:= \pi_{r-1}(M)$).
Hence, $E2=E^r$, and $E^{r+1}=E^{\infty}$.

Since the spectral sequence converges to $H_{s+t}(M; R)$,
and $H_r(M; R)=0$, we infer that
\begin{equation}
\label{eq:dmono}
d^r : E_{r0}^r =\tor_r^{\Z \pi}(\Z, R)\hookrightarrow
\pi_{r-1}\otimes_{\nu} R= E_{0, r-1}^r
\end{equation}
is monic. Convergence also provides the exact sequence
\begin{equation}
\label{eq:conv}
0\to E_{0, r-1}^{\infty}=\coker \, \{ d^r \} \longrightarrow
H_{r-1} (M; R) \longrightarrow E_{r-1, 0}^{\infty} =
\tor_{r-1}^{\Z \pi} (\Z, R)\to 0 \, .
\end{equation}

Looking at $R$--ranks in~\eqref{eq:conv}, and also
taking into account~\eqref{eq:dmono}, we find that
\begin{equation}
\label{eq:fin}
\rank_R \pi_{r-1}\otimes_{\nu} R=
\rank_R H_{r-1}(M; R)-
\rank_R \tor_{r-1}^{\Z \pi}(\Z ,R)+
\rank_R \tor_{r}^{\Z \pi}(\Z ,R) \, .
\end{equation}

We may now substitute in equation~\eqref{eq:fin} above 
$\rank_R H_{r-1}(M(\A); R)$
by the formula from Theorem~\ref{thm:compl}\eqref{c1},
to get the asserted formula for
$\rank_R \pi_{r-1}(M(\A))\otimes_{\nu} R$.
\end{proof}

The $R$--modules, $\pi_{r-1}(M(\A))\otimes_{\nu} R$
from Theorem~\ref{thm:rarb}, and 
$ H_{r-1}(M(\A); R)$ from Theorem~\ref{thm:compl}\eqref{c1},
look similar. Indeed,
$\pi_{r-1}(M(\A))\otimes_{\nu} R = H_{r-1}(C_{\bullet})
\otimes_{\Z \pi} R$, and
$ H_{r-1}(M(\A); R)= H_{r-1}(C_{\bullet}\otimes_{\Z \pi} R)$,
where
$C_{\bullet}:= C_{\bullet}(\tY^{(r-1)})$, and
$Y$ is a minimal structure on $M(\wA)$. In fact, they
turn out to be different, in general, even at the level of
$R$--ranks.

\begin{example}
\label{ex:braid}

Let $\wA$ be the essential arrangement in $\PP^r$ associated to
the braid arrangement from $\PP^{r+1}$, $r\ge 3$. Let $\A$ be a
generic (i.e., $\LL_{r-1}(\wA)$--generic) hyperplane section 
of $\wA$. Let $\pi$ be the pure braid group on $r+2$ strings
($\pi = \pi_1(M(\A))= \pi_1(M(\wA))$). Let $\nu : \pi\to \Z$ be
$\wA$--nonresonant, giving rise to the local system 
$R:= \C \Z$.

Then $\rank_R \pi_{r-1}(M(\A))\otimes_{\nu} R=
\rank_R  H_{r-1}(M(\A); R)$ if and only if
\begin{equation}
\label{eq:top2}
\rank_R H_{r-1}(Y; R)=
\rank_R H_r(Y; R)\, ,
\end{equation}
where $Y$ is a minimal structure on $M(\wA)$; see
Theorem~\ref{thm:rarb} and Theorem~\ref{thm:compl}\eqref{c1}.

By Remark~\ref{rem:tors}, \eqref{eq:top2} holds if
and only if $\chi (M(\wA))=0$. This in turn is impossible,
since the well-known Poincar\' e polynomial of $M(\wA)$ is
$\prod_{d=2}^{r+1} (1+ dT)$.
\end{example}

\noindent {\bf Acknowledgement.} This work was completed while
the second author was visiting the University of Bordeaux. He
thanks the CNRS and the Bordeaux Department of Mathematics for
the excellent working facilities provided.


\begin{thebibliography}{00}

\bibitem{C} D. Cohen:
Cohomology and intersection cohomology of complex
hyperplane arrangements,
{\em Adv. Math.} \textbf{97} (1993), 231--266.

\bibitem{CDO} D. Cohen, A. Dimca and P. Orlik:
Nonresonance conditions for arrangements,
{\em preprint} {\tt math.AG/0210409}.

\bibitem{CO} D. Cohen and P. Orlik:
Arrangements and local systems,
{\em Math. Research Letters} \textbf{7} (2000),
299--316.

\bibitem{CS1} D. Cohen and A. Suciu:
On Milnor fibrations of arrangements,
{\em J. London Math. Soc.} \textbf{51} (1995),
105--119.

\bibitem{CS3} D. Cohen and A. Suciu:
The braid monodromy of plane algebraic curves 
and hyperplane arrangements,
{\em Comment. Math. Helv.} \textbf{72} (1997), 285--315.

\bibitem{CS2} D. Cohen and A. Suciu:
Homology of iterated semidirect products of free groups,
{\em J. Pure Appl. Algebra} \textbf{126} (1998), 87--120.

\bibitem{DN} A. Dimca and A. N\' emethi:
Hypersurface complements, Alexander modules and monodromy,
{\em preprint} {\tt math.AG/0201291}.

\bibitem{DP} A. Dimca and S. Papadima:
Hypersurface complements, Milnor fibers and minimality
of arrangements,
{\em preprint} {\tt math.AT/0101246},
to appear in 
{\em Ann. Math.} \textbf{157} (2003).

\bibitem{E} D. Eisenbud:
``Commutative algebra with a view toward algebraic geometry'',
{\em Grad. Texts in Math.}, Vol. \textbf{150},
Springer-Verlag, New York, 1995.

\bibitem{FR} M. Falk and R. Randell:
The lower central series of a fiber-type arrangement,
{\em Invent. Math.} \textbf{82} (1985), 77--88.

\bibitem{GWPL}  C. G. Gibson, K. Wirthm\" uller, 
A. A. du Plessis and E. J. N. Looijenga: 
``Topological stability of smooth mappings'', 
{\em LNM} \textbf{552}, Springer-Verlag, Berlin, 1976.

\bibitem{GM} M. Goresky and R. MacPherson:
``Stratified Morse theory'',
{\em Ergebnisse}, Vol. \textbf{14},
Springer-Verlag, New York, 1988.

\bibitem{H} A. Hattori:
Topology of $\C^n$ minus a finite number of affine hyperplanes 
in general position,
{\em J. Fac. Sci. Univ. Tokyo} \textbf{22} (1975), 205--219.

\bibitem{Hi} J. A. Hillman:
``Alexander ideals of links'', {\em LNM} \textbf{895},
Springer-Verlag, Berlin, 1981.

\bibitem{JP1} M. Jambu and S. Papadima:
A generalization of fiber-type arrangements and a new deformation
method,
{\em Topology} \textbf{37} (1998), 1135--1164.

\bibitem{JP2} M. Jambu and S. Papadima:
Deformations of hypersolvable arrangements,
{\em Topology Appl.} \textbf{118} (2002), 103--111.

\bibitem{L1} A. Libgober:
On the homotopy type of the complement to plane
algebraic curves, {\em J. Reine Angew. Math.}
\textbf{397} (1986), 103--114.

\bibitem{L2} A. Libgober:
Homotopy groups of the complements to singular hypersurfaces II,
{\em Ann. Math.} \textbf{139} (1994), 117--144.

\bibitem{ML} S. Mac Lane: 
``Homology'',
{\em Grundlehren}, Vol. \textbf{114},
Springer-Verlag, Berlin, 1963.

\bibitem{OS} P. Orlik and L. Solomon:
Combinatorics and topology of complements of hyperplanes,
{\em Invent. Math.} \textbf{56} (1980), 167--189.

\bibitem{OT} P. Orlik and H. Terao:
``Arrangements of hyperplanes'',
{\em Grundlehren}, Vol. \textbf{300},
Springer-Verlag, Berlin, 1992.

\bibitem{PS} S. Papadima and A. Suciu:
Higher homotopy groups of complements of complex 
hyperplane arrangements,
{\em Adv. Math.} \textbf{165} (2002), 71--100.

\bibitem{R} R. Randell:
Morse theory, Milnor fibers and minimality of
hyperplane arrangements, {\em Proc. Amer. Math. Soc.} 
\textbf{130} (2002), 2737--2743.

\bibitem{Ry} G. Rybnikov:
On the fundamental group of the complement of a complex
hyperplane arrangement, {\em DIMACS Tech. Report}
\textbf{94-13} (1994), 33--50; available at
{\tt math.AG/9805056}.
 
\bibitem{STV} V. Schechtman, H. Terao and A. Varchenko:
Local systems over complements of hyperplanes and the
Kac--Kazhdan condition for singular vectors,
{\em J. Pure Appl. Alg.} \textbf{100} (1995), 93--102.

\bibitem{W} G. W. Whitehead:
``Elements of homotopy theory'',
{\em Grad. Texts in Math.}, Vol. \textbf{61},
Springer-Verlag, New York, 1978.
 
\end{thebibliography}
\end{document}